\def\draft{n}

\documentclass[11pt]{amsart}

\usepackage{amssymb, amsmath, latexsym, amscd}
\usepackage{graphicx}
\usepackage{mathpple}

\def\printname#1{
	\if\draft y
		\smash{\makebox[0pt]{\hspace{-0.5in}
			\raisebox{8pt}{\tt\tiny #1}}}
	\fi }
\def\lbl#1{\label{#1}\printname{#1}}

\newtheorem{theorem}{Theorem}[section]
\newtheorem{corollary}[theorem]{Corollary}

\newtheorem{conjecture}[theorem]{Conjecture}
\newtheorem{thm}[theorem]{Theorem}

\newtheorem{prop}[theorem]{Proposition}

\theoremstyle{definition}

\begin{document}

\newcommand{\Ass}{\mathcal A}ÊÊÊÊÊÊÊÊÊÊÊÊÊ 
\newcommand{\Com}{ \mathcal C}ÊÊÊÊÊÊÊÊÊÊ 
\newcommand{\F}{\mathcal F}
\newcommand{\Lie}{\mathcal L}ÊÊÊÊÊÊÊÊÊÊÊÊÊÊ
\newcommand{\bdry}{\partial}
\newcommand{\g}{\mathcal G}
\newcommand{\fox}{d}
\newcommand{\assfox}{d^{\mathbf a}}
\newcommand{\liefox}{\mathfrak D}
\newcommand{\rpartial}{D}
\newcommand{\ls}{\mathcal{LL}} 
\newcommand{\la}{\mathcal{LA}}
\newcommand{\lc}{\mathcal{LC}} 
\newcommand{\tr}{\Tr^{\mathcal G}}
\newcommand{\Tr}{\operatorname{Tr}}
\renewcommand{\b}{\mathbf B}
\newcommand{\B}{\mathbb B}
\newcommand{\Q}{\mathbb Q}
\newcommand{\m}{\mathcal G_{AB}}
\renewcommand{\o}{o}
\newcommand{\Ge}{\mathbf G_e}
\newcommand{\GC}{\mathbf G\langle C\rangle}
\newcommand{\zexp}{\operatorname{exp}_{\mathbb Z}}
\newcommand{\zlog}{\operatorname{log}_{\mathbb Z}}

\title[Chirality conjecture]{Chirality and the Conway polynomial}
\author{James Conant}
\address{Department of Mathematics\\
ÊÊÊÊÊÊÊÊÊ University of Tennessee\\
ÊÊÊÊÊÊÊÊÊ Knoxville, TN, 37996}
\email{jconant@math.utk.edu}

\begin{abstract} 
In recent work with J. Mostovoy and T. Stanford, the author found that for every natural number n, a certain polynomial in the coefficients of the Conway polynomial is a primitive integer-valued degree n Vassiliev invariant, but that modulo 2, it becomes degree n-1.
The conjecture then naturally suggests itself that these primitive invariants are congruent to integer-valued degree n-1 invariants. In this note, the consequences of this conjecture are explored. Under an additional assumption, it is shown that this conjecture implies that 
the Conway polynomial of an amphicheiral knot has the property that $C(z)C(iz)C(z^2)$ is a perfect square inside the ring of power series with integer coefficients, or, equivalently, the image of $C(z)C(iz)C(z^2)$ is a perfect square inside the ring of polynomials with $\mathbb Z_4$ coefficients.  In fact, it is probably the case that the Conway polynomial of an amphicheiral knot always can be written as $f(z)f(-z)$ for some polynomial $f(z)$ with integer coefficients, and this actually implies the above ``perfect squares" conditions. Indeed, by work of Kawauchi and Hartley, this is known for all negative amphicheiral knots and for all strongly positive amphicheiral knots. In general it remains unsolved, and this paper can be seen as
some evidence that it is indeed true in general.
\end{abstract}

\maketitle

\section{Introduction}
The Conway polynomial of knots, which is a renormalization of the Alexander polynomial, is of interest partially because its coefficients $c_n$ are Goussarov-Vassiliev invariants of degree $n$. The invariants $c_n$ are not, however, additive under the connected sum operation, which is inconvenient, because the analysis of weight systems is considerably easier in the additive, or \emph{primitive}, case. 
Since $C(z)$ is multiplicative under the connected sum operation,
one can apply a standard trick and take the formal logarithm of $C(z)$ regarded as a power series, but this turns integer coefficients into rational ones, making the study of torsion harder. In \cite{CMS}, we obviated this problem by introducing a discrete logarithm $\zlog\colon 1+x\cdot\mathbb Z[[x]]\to x\cdot\mathbb Z[[x]]$ which is a bijection sending multiplication to addition. If we write
$$\zlog C(z)=\sum _{n=1}^\infty pc_{2n}z^{2n}$$
then $pc_{2n}$ are primitive Goussarov-Vassiliev invariants of degree $2n$, and moreover, modulo $2$, they are of degree $2n-1$, except in the case $n=1$. 

This is one piece of evidence to support the conjecture (Conjecture \ref{gv}) that the invariants $pc_{4n}$ are congruent modulo $2$ to degree $4n-1$ invariants $v_{4n-1}$. 
(The conjecture is stated only for $4n$, although
it may be that $pc_{2n}$ ($n\geq 2$) are also congruent to some invariants $v_{2n-1}$.) 
 Conjecture~\ref{gv} follows from the slightly stronger Conjecture~\ref{strange2} which asserts that the mystery invariants $v_{4n-1}$ change sign when a knot is sent to its mirror image.  (For example, the odd degree part of the Kontsevich integral has this property.) This stronger conjecture has the following surprising consequence.

\begin{conjecture}\lbl{main}
Suppose $K$ is an amphicheiral knot. Then there is a polynomial $F\in\mathbb Z_4[z^2]$ such that
$$F^2=C(z)C(z^2)C(iz)\in\mathbb Z_4[z^2].$$
\end{conjecture}

This conjecture holds for a wide class of amphicheiral knots including all amphicheiral knots in the knot tables with 14 and fewer crossings. (There are several hundred such knots.)
In fact, more seems to be true. Let us temporarily adopt the terminology that the Conway polynomial is \emph{splittable} if it can be factored  as $f(z)f(-z)$ for some polynomial $f(z)$ with integer coefficients. Then
Proposition~\ref{HKProp} indicates that if the Conway polynomial of an amphicheiral knot is always splittable, then Conjecture~\ref{main} is true.
Hartley has shown \cite{hartley} that every negative amphicheiral knot has a splittable Conway polynomial. (A knot is negative amphicheiral if it is isotopic to its string-orientation reversed mirror image.) This built on earlier work of Hartley and Kawauchi \cite{HK} which established this for strongly amphicheiral knots. (A knot is strongly amphicheiral if it is fixed set-wise by an orientation-reversing involution of $S^3$.)
The general case, the splittability of the Conway polynomial of an arbitrary amphicheiral knot, remains unsolved. Besides the results of this paper,
there is another bit of evidence that it is true. Namely, there is a result of L. Goeritz \cite{G} which implies that the determinant of an amphicheiral knot is the sum of two squares. As pointed out by Hartley and Kawauchi, this is implied by the fact that 
$C(z)=f(z)f(-z)$. Since the determinant can be calculated by setting $z=2i$, the result is a complex number times its conjugate: a sum of two squares. 

It is not hard to show, using Mostow rigidity, that any hyperbolic amphicheiral knot which has $\mathbb Z_{4k}$ as a symmetry group cannot be negative amphicheiral nor can it be strongly amphicheiral, so that this is a source of possible counterexamples to the splittability conjecture. However, all such knots the author has tested have satisfied $C(z)=f(z)f(-z)$. In fact, these examples are in some sense atomic. Hartley proves that all negative amphicheiral knots have splittable Conway polynomial by considering the JSJ decomposition of the knot complement.
The argument proceeds inductively, with the base of the induction being the case of knots with complements which are either Seifert fibred of hyperbolic. There are no amphicheiral seifert fibred knots, so we are left with hyperbolic knots. It is well known that the symmetry group of a hyperbolic knot is either cyclic or dihedral. Using Mostow rigidity, the amphicheiral symmetry must be realizable by a hyperbolic symmetry. Then in all cases but $\mathbb Z_{4k}$ it is easy to produce an orientation reversing involution, implying the knot is strongly amphicheiral. Then one appeals to the main result of \cite{HK}. The cases of $\mathbb Z_{4k}$ do not arise in Hartley's work since such knots are not negative amphicheiral.

{\bf Acknowledgments:}
I wish to thank I. Agol, R. Budney, C. Livingston, J. Mostovoy, S. Mulay, T. Stanford and M. Thistlethwaite for contributing to the ideas in this paper in one way or another. The author was partially supported by NSF grant DMS 0305012.

\section{A Formal Logarithm}
We begin by defining a formal exponential
$$\zexp\colon x\cdot \mathbb Z[[x]]\to 1+x\cdot \mathbb Z[[x]],$$
where $\zexp(F+ G)= \zexp(F)\cdot\zexp(G)$.
Namely, let
$$\zexp\left(\sum_{i=1}^\infty a_i x^i\right)=\prod_{i=1}^\infty(1+(-x)^i)^{a_i}.$$
Evidently
$\zexp$ takes addition to multiplication. 

It is also rather easy to see that $\zexp$ is a bijection. This follows readily from the fact that $$\zexp\left(\sum_{i=1}^\infty a_i x^i\right)=1+\sum_{i=1}^\infty
(a_i+p_i(a_1,\ldots,a_{i-1}))x^i,$$
where $p_i$ is a polynomial with integer coefficients.

Now define $$\zlog\colon\! 1+x\cdot \mathbb Z[[x]]\to x\cdot \mathbb Z[[x]]$$ as the inverse to $\zexp$. By construction, it takes multiplication to addition.

\begin{proof}[Example:] $$\zlog(1+x)=\sum_{i=1}^\infty -x^{2^i}.$$This is equivalent to saying that
$$1+x=\prod_{i=0}^\infty(1+(-x)^{2^i})^{-1},$$ which can be proven using an elementary telescoping argument. 
\end{proof}

If we write $\zlog \left(1+\sum_{i=1}^\infty a_i x^i\right)=\sum_{i=0}^\infty b_i x^i$, then formulae for the first few $b_i$ are given as follows
\begin{align*}
b_1&= -a_1\\
b_2&= a_2-\frac{1}{2}(a_1+a_1^2) \\
b_3&= -a_3+a_1a_2+\frac{1}{3}(a_1-a_1^3)\\
b_4&= a_4-a_1a_3+\frac{1}{2}(a_2-a_2^2)+a_1^2a_2-\frac{1}{4}(2a_1+a_1^2+a_1^4)
\end{align*}

\section{Primitive Goussarov-Vassiliev invariants}

The algebra of rational-valued Goussarov-Vassiliev invariants is actually a Hopf algebra, and so by the structure theory of Hopf algebras, it is a polynomial algebra generated by primitive elements. The condition of primitivity amounts to saying that a knot invariant is additive under connected sum. Thus primitive knot invariants play a special role in the theory, although it's not clear whether they polynomially generate all Goussarov-Vassiliev invariants in the case where the invariant's target group is more general. 

The coefficients of the Conway polynomial $c_{2i}$ are Goussarov-Vassiliev invariants of degree $2i$, but they are not primitive.
By applying $\zlog$ to $C(z^2)$ (regarding $z^2$ as the variable) we get an additive invariant of knots taking values in
$\mathbb Z[[z^2]]$, and it is easy to show that the coefficients, denoted $pc_{2i}$, are primitive Vassiliev invariants of degree $2i$. 

The following theorem is proven in \cite{CMS}.
\begin{thm}\lbl{strange}
$pc_{2i}$ is a degree $2i$ Goussarov-Vassiliev invariant over $\mathbb Z$, and is of degree
$2i-1$ over $\mathbb Z_2$.
\end{thm}

In that paper we deduce the consequence that the only primitive Goussarov-Vassiliev invariants of $S$-equivalence come from the invariants $pc_{2i}$. (Two knots are $S$-equivalent if they have isomorphic Alexander modules and Blanchfield forms. Alternatively, if there exist Seifert surfaces with isomorphic Seifert pairings.)

Theorem \ref{strange} opens the door for us, because it suggests that there may be an integer-valued degree $2i-1$ invariant $v_{2i-1}$ which is a lift of $pc_{2i}$ modulo $2$. There isn't enough evidence to conjecture that all of these $pc_{2i}$ lift, however, there is some evidence to support the following.

\begin{conjecture}\lbl{gv}
There exist integer-valued Goussarov-Vassiliev invariants, $v_{4i-1}$, of degree $4i-1$,
such that
$$pc_{4i}\equiv v_{4i-1}\mod 2.$$
\end{conjecture}
 
Something that often happens is that odd-degree invariants are odd. That is, they change sign under mirror image, and therefore vanish on amphicheiral knots. (Conceivably string orientation
might be a factor, but no known Goussarov-Vassiliev invariant can detect string orientation. Moreover the fact that the amphicheirality criterion we ultimately derive works for $+/-$ amphicheiral knots
 suggests that the invariants $v_{4i-1}$, if they exist, are probably insensitive to string orientation. )

Therefore, if $v_{2i-1}$ were indeed an odd invariant, then $pc_{2i}$ would have to be even on any amphicheiral knot. Computer experiment quickly disposes of this possibility for odd $i$. However, computer evidence does support the following conjecture.

\begin{conjecture}\lbl{strange2}
The invariants $v_{4i-1}$ of Conjecture~\ref{gv} can be chosen to be odd. 
\end{conjecture} 

This conjecture is a fact proven by Ted Stanford when $i=1$ \cite{S}, where it is shown that a degree $3$ odd invariant, $v_3$, is
congruent, modulo 2, to $\frac{1}{2}(c_2+c_2^2)+c_4.$ Indeed, we may take
$v_3=-\frac{1}{12}J''(1)-\frac{1}{36}J'''(1)$, where $J(t)$ is the Jones polynomial. 

\begin{corollary}\lbl{crit}
If Conjecture~\ref{strange2} holds, then  $pc_{4i}\equiv 0\mod 2$ on amphicheiral knots.
\end{corollary} 

This gives us a sequence of amphicheirality criteria. We list the first three in terms of the standard coefficients of the Conway polynomial.

\begin{align*}
&c_4+\frac{1}{2}(c_2^2+c_2)\equiv 0\mod 2\text{\hspace{13em}{(Stanford's Criterion)}}\\
&c_8+c_2(c_6+c_4)+\frac{1}{2}(c_4^2-c_4)+\frac{1}{4}(c_2^4+c_2^2+2c_2)\equiv 0\mod 2\\
&c_{12}+c_2(c_{10}+c_4+c_6+c_8)+c_4c_8+\frac{1}{2}(c_6^2+c_6)+
\frac{1}{2}(c_2c_4-3(c_2c_4)^2)\equiv 0\mod 2
\end{align*}

 \section{A conjecture without a formal log}
Conjecture \ref{strange} is certainly intriguing and the amphicheirality criterion of Corollary~\ref{crit} can be easily tested with the aid of a computer mathematics package. However, it turns out that the criterion can be reformulated without reference to $\zlog$. This is shown in the following theorem.
 
 \begin{thm}\lbl{stuff}
Let $$f(x)=\prod_{i=1}^\infty (1+(-x)^i)^{a_i}.$$
Then $f(x)f(-x)f(x^2)\in (1+x\cdot\mathbb Z[[x]])^2$ if and only if for all $i$, $a_{2i}\equiv 0\mod 2$.
 \end{thm}
\begin{proof}
Suppose $f(x)=\prod_{i=1}^\infty (1+x^i)^{a_i}$ satisfies $a_{2i}\equiv 0$ for all $i$. Then  we wish to show that $f(x)f(-x)f(x^2)$ is a square inside
$1+x\cdot \mathbb Z[[x]]$. Then $f(x)$ is equivalent modulo squares to 
$\prod_{i=1}^\infty(1-x^{2i-1})^{a_{2i-1}}.$ Now 
$$f(x)f(-x)f(x^2)\equiv\prod_{i=1}^\infty (1-x^{2i-1})^{a_{2i-1}}(1+x^{2i-1})^{a_{2i-1}}(1+x^{4i-2})^{a_{2i-1}},$$
which is a product of squares.

For the converse, note that the previous argument shows that
\begin{align*}
f(x)f(-x)f(x^2)&\equiv\prod_{i=1}^\infty (1+x^{2i})^{a_{2i}}(1+x^{2i})^{a_{2i}}(1+x^{4i})^{a_{2i}}\\
&\equiv\prod_{i=1}^\infty (1+x^{4i})^{a_{2i}}
\end{align*}
Suppose, toward a contradiction that not all of the $a_{2i}$ are even.
Let $i_0$ be the minimal index such that $a_{2i_0}$ is odd. Then modulo squares,
\begin{align*}
f(x)f(-x)f(x^2)&\equiv \prod_{i=i_0}^\infty(1+x^{4i})^{a_{2i}}\\
&= 1+a_{2i_0}x^{4i_0}+\cdots
\end{align*}
Since the coefficient of the smallest power of $x$ must be even in a perfect square, we deduce that $a_{2i_0}$ is even, a contradiction.
\end{proof} 

The next two propositions allow us to deduce Conjecture~\ref{main} from Conjecture~\ref{strange2} and Theorem~\ref{stuff}.

  \begin{prop}\lbl{stuff2}
 Suppose $F, G\in1+x\cdot\mathbb Z[[x]]$, and the coefficients of $F$ are congruent to the coefficients of $G$ modulo $4$.
 If $\sqrt{F}\in 1+x\cdot\mathbb Z[[x]],$ then $\sqrt{G}\in1+x\cdot\mathbb Z[[x]]$.
 \end{prop}
\begin{proof}
Assume that $G=F\pm 4x^i$. This will imply the general case. 
I claim that $\sqrt{G/F}\in 1+x\cdot\mathbb Z[[x]]$ in which case we're done, since
$\sqrt{G}=\sqrt{G/F}\sqrt{F}$.  Letting $\gamma=\pm x^i/F$, we have that
$$\sqrt{G/F}=\sqrt{1+4\gamma}$$
It is now an exercise in Taylor series to show that $\sqrt{1+4\gamma}$ has integer coefficients.
\end{proof}

\begin{prop}\lbl{stuff3}
Let $F$ be a formal power series in $\mathbb Z_4[[x]]$. Suppose $F^2\in\mathbb Z_4[x]$. That is, the square is a polynomial. Then either $F\in \mathbb Z_4[x]$ or $F=F_1+2F_2$ where $F_1$ is a polynomial.
\end{prop}
\begin{proof}
Suppose otherwise. Then $F$ has coeffcients of $\pm 1$ of arbitrarily high degree.
Consider $\bar F\in \mathbb Z_2[[x]]]$. Then $\bar F=\sum_{i\in I} z^i$, where $|I|$ is infinite. Now $\bar F^2=\sum_{i\in I} z^{2i}$ which is not a polynomial, indicating that $F^2$ is also not a polynomial.
\end{proof}

\begin{prop}
Conjecture~\ref{main} is equivalent to the criterion of Corollary~\ref{crit}.
\end{prop}
\begin{proof}
Suppose $pc_{2i}\equiv pc_{4i}$ for all $i$. According to Theorem ~\ref{stuff}, this is equivalent to 
$C(z)C(iz)C(z^2)$ being the square of a formal power series with integer coefficients. By Proposition~\ref{stuff2}, this is equivalent to the corresponding statement with $\mathbb Z_4$ coefficients.  Now, if $C(z)C(z^2)C(iz)=F^2$ inside $\mathbb Z_4[[z]]$, then we know $F=F_1+2F_2$ by Proposition \ref{stuff3}, where $F_1$ is a polynomial. But then $C(z)C(z^2)C(iz)=F_1^2$.
\end{proof}

\section{The work of Hartley and Kawauchi}
Recall that a knot is said to be strongly positive amphicheiral if there is an orientation reversing involution of $S^3$ which fixes the knot set-wise, and whose restriction to the knot preserves string orientation.

Hartley shows that if a knot is negative amphicheiral, then its Alexander polynomial can be written $f(\sqrt{t})f(-\sqrt{t})$ for some integral polynomial $f(\sqrt{t})$ satisfying
$f(-\sqrt{t})=f(1/\sqrt{t})$ and $|f(1)|=1$. This can be rephrased as the fact that
the Conway polynomial factors as $C(z)=\phi(z)\phi(-z)$, where $\phi$ is an integer polynomial.

Similarly, Hartley and Kawauchi show that if a knot is strongly positive amphicheiral, then its Alexander polynomial can be written $f(t)^2$ for some integral polynomial $f(t)$ satisfying
$f(t)=f(t^{-1})$ and $|f(1)|=1$. This can be rephrased as the fact that the Conway polynomial factors as
$C(z)=\psi(z^2)^2$ for an integer polynomial $\psi.$ In particular, the Conway polynomial is of the above form where $\phi(z)=\psi(z^2)$.

\begin{prop}\lbl{HKProp}
If $C(z)=\phi(z)\phi(-z)$ for an integer polynomial $\phi$, then $C(z)C(iz)C(z^2)$ is the square of an integer power series. Hence Conjecture~\ref{main} holds for all strongly amphicheiral knots.
\end{prop}
\begin{proof}
Write $C(z)=\phi(z)\phi(-z)$. Either $\phi(0)=1$ or $\phi(0)=-1$. If the latter, replace $\phi$ by $-\phi$.
Then $\phi(z)=\prod_{i=1}^\infty (1+(-z)^i)^{a_i}.$ So, modulo squares, $\phi(z)\phi(-z)= \prod_{i=1}^\infty (1-z^{2i-1})^{a_{2i-1}}$. On the other hand $C(iz)$ is modulo squares the product of terms of the form $(1+z^{2i-1})^{a_{2i-1}}$, and $C(z^2)$ is modulo squares a product of terms of the form $(1-z^{4i-2})^{a_{2i-1}}$, from which the desried result follows easily.
 \end{proof}

The condition that $C(z)=\phi(z)\phi(-z)$ is strictly stronger than that of Conjecture ~\ref{main}, even if one assumes the determinant is a sum of two squares. For example $C(z)=1-76z^2$ is congruent to $1$ modulo $4$, so satisfies the criterion of Conjecture~\ref{main}, and $\operatorname{det}=C(2i)=305$ is a sum of two squares since the primes in its factorization are congruent to $1$ modulo $4$. However $C(z)\neq\phi(z)\phi(-z)$ since the absolute value of the coefficient of the highest power of $z$ is not a square.

\end{document}